\documentclass[runningheads]{amsart}

\usepackage[utf8]{inputenc}

\usepackage{amssymb}

\usepackage[T1]{fontenc}

\newcommand{\Ld}[2][D]{\ensuremath{\mathcal{L}(#2\mathcal{#1})}} 
\newtheorem{theorem}{Theorem}[section]
\newtheorem{lemma}[theorem]{Lemma}
\newtheorem{proposition}[theorem]{Proposition}
\newtheorem{corollary}[theorem]{Corollary}
\theoremstyle{definition}
\newtheorem{definition}[theorem]{Definition}
\newtheorem{remark}[theorem]{Remark}

\newcommand{\bbbf}{\mathbb{F}}
\newcommand{\bbbn}{\mathbb{N}}

\begin{document}
\title[Chaining multiplications]{Chaining Multiplications
in Finite Fields with Chudnovsky-type Algorithms and Tensor Rank of the k-multiplication}
\author{S. Ballet \and R. Rolland}
\address{Aix Marseille Univ, CNRS, Centrale Marseille, I2M, Marseille, France \\ Case 907
13288 MARSEILLE Cedex 9\\ France}

\maketitle

\begin{abstract}
We design a class of Chudnovsky-type algorithms multiplying $k$ elements
of a finite extension ${\bbbf}_{q^n}$ of a finite field ${\bbbf}_q$, where $k\geq 2$.
We prove that these algorithms give a tensor decomposition of the k-multiplication for which the rank is in $O(n)$ uniformly in $q$. 
We give uniform upper bounds of the rank of k-multiplication in finite fields.
They use interpolation on algebraic curves which transforms the problem
in computing the Hadamard product of $k$
vectors with components in ${\bbbf}_q$. 
This generalization of
the widely studied case of $k=2$
is based on a modification 
of the Riemann-Roch spaces involved and the use of towers of function fields having a lot of places of high degree.

\keywords{Finite Field \and Tensor Rank of the Multiplication \and Function Field}
\end{abstract}

%\tableofcontents

\section{Introduction}

\subsection{Context, notation and basic results}

The construction  of efficient arithmetic operation algorithms is still a problem of topicality. 
These algorithms are indeed heavily used in many domains of computer sciences or information theory.  
 It is important to conceive and develop efficient arithmetic algorithms combined with an optimal implementation method. 
 In this work, our interest lies in multiplication algorithms in any extension of finite field introduced in 1987 by  D.V. and G.V Chudnovsky \cite{chch} 
 and based upon interpolation on points of algebraic curves defined over finite fields. 
 Our goal is to improve this method so that its complexity in terms of number of operations be optimized. 
 
 \medskip

More precisely, the complexity of a multiplication algorithm in ${\bbbf}_{q^n}$ depends on the number
of multiplications and additions in ${\bbbf}_q$. But here, we are particularly interested by the multiplicative complexity 
of multiplication in a finite field ${\bbbf}_{q^n}$, i.e. by the number of multiplications in ${\bbbf}_q$ required to multiply in the 
${\bbbf}_q$-vector space ${\bbbf}_{q^n}$ of dimension $n$. There exist two types of multiplications in ${\bbbf}_q$: 
the scalar multiplication and the  bilinear  one. The scalar multiplication is the multiplication by a non-trivial constant 
(i.e. not equal to $0$ or $1$) in ${\bbbf}_{q}$, which does not depend on the elements of ${\bbbf}_{q^n}$
that are multiplied. The bilinear multiplication is a multiplication that depends on the elements of ${\bbbf}_{q^n}$ 
that are multiplied. The bilinear complexity is independent of the  chosen representation of the finite field.

 \medskip

Let $q$ be a prime power, ${\bbbf}_q$ the finite field with $q$ elements and ${\bbbf}_{q^n}$ the degree $n$ extension of ${\bbbf}_q$.
If $\mathcal{B}=\{e_1,...,e_n\}$ is a basis of ${\bbbf}_{q^n}$ over ${\bbbf}_q$ then for $x=\sum_{i=1}^{n}x_ie_i$ and $y=\sum_{j=1}^{n}y_je_j$, we have the product
\begin{equation}\label{calculdirect}
z=xy=\sum_{h=1}^{n}z_he_h=\sum_{h=1}^{n}\biggr( \sum_{i,j=1}^{n}t_{ijh}x_iy_j\biggl)e_h,
\end{equation}
 where $e_ie_j=\sum_{h=1}^{n}t_{ijh}e_h,$ $t_{ijh}\in {\bbbf}_q$ being some constants.
 
\vspace{1em}

Then, we see that the direct calculation of $z=(z_1,...,z_n)$ using (\ref{calculdirect}) {\it a priori} requires $n^2$ non-scalar multiplications $x_iy_j$, $n^3$  scalar multiplications and $n^3-n$ additions.

\begin{definition}\label{defmus}
The total number of scalar multiplications in ${\bbbf}_q$ used in an algorithm ${\mathcal U}_{q,n}$ of multiplication in ${\bbbf}_{q^n}$ is
called scalar complexity of ${\mathcal U}_{q,n}$ and denoted $\mu_S({\mathcal U}_{q,n})$.
\end{definition}

\vspace{1em}

Moreover, 
if $k$ is an integer $\geq 2$, the multiplication $m_k$ of $k$ elements in the finite field ${\bbbf}_{q^n}$ is a $k$-multilinear map 
from $\left({\bbbf}_{q^n}\right)^k$ into ${\bbbf}_{q^n}$ over the field ${\bbbf}_q$,
thus it corresponds to a linear map $M_k$ from the tensor power
$({\bbbf}_{q^n})^{\otimes k}$ into ${\bbbf}_{q^n}$.
One can also represent $M_k$ by a $k$-covariant and $1$-contravariant tensor 
$t_{M_k} \in ({\bbbf}_{q^n}^*)^{\otimes k} \otimes {\bbbf}_{q^n}$ 
where ${\bbbf}_{q^n}^*$ denotes the algebraic dual of ${\bbbf}_{q^n}$.
Each decomposition
\begin{equation}\label{algo}
t_{M_k}=\sum_{i=1}^{s} \left(\otimes_{j=1}^k a^*_{i,j}\right) \otimes c_i
\end{equation}
of the tensor $t_{M_k}$, where $a^*_{i,j} \in {\bbbf}_{q^n}^*$ and
$c_i \in {\bbbf}_{q^n}$, brings forth a multiplication algorithm of $k$ elements
\begin{equation}\label{k-multilinearAlgo}
\prod_{j=1}^k x_j=t_{M_k}(\otimes_{i=1}^k x_j)=\sum_{i=1}^{s}\left(\prod_{j=1}^ka^*_{i,j}(x_j)\right) c_i .
\end{equation}

\begin{definition}
A $k$-multilinear multiplication algorithm ${\mathcal U}_{q,n,k}$ in ${\bbbf}_{q^n}$
is an expression

$$\prod_{j=1}^k x_j=\sum_{i=1}^{s}\left(\prod_{j=1}^ka^*_{i,j}(x_j)\right) c_i .$$
where $a^\star_{i,j} \in ({{\bbbf}{q^n}})^\star$, and $c_i\in {\bbbf}{q^n}$.

The number $s$ of summands in this expression is called the $k$-multilinear complexity of the algorithm ${\mathcal U}_{q,n,k}$
and is denoted by $\mu_M({\mathcal U}_{q,n,k})$. The multiplicative complexity of ${\mathcal U}_{q,n,k}$ is $\mu_{m}({\mathcal U}_{q,n,k})= \mu_M({\mathcal U}_{q,n,k})+\mu_S({\mathcal U}_{q,n,k})$.

\end{definition}

\begin{definition}
The $k$-multilinear complexity of the multiplication of $k$ elements in ${\bbbf}_{q^n}$ over ${\bbbf}_q$, 
which is denoted by $\mu_{q,k}(n)$,
is the minimum number $s$ of summands in the decomposition (\ref{algo}). 
\end{definition}

\begin{definition}
The minimal number of summands in a decomposition of the tensor $T_{M_k}$ of the $k$-multilinear multiplication in ${\bbbf}_{q^n}$
is called the $k$-multilinear complexity of the multiplication in ${\bbbf}_{q^n}$ and is denoted by
$\mu_{q,k}(n)$: 
$$
\mu_{q,k}(n)= \min_{{\mathcal U}_{q,n,k}} \mu_M({\mathcal U}_{q,n,k}) 
$$
where ${\mathcal U}_{q,n,k}$ is running over all $k$-multilinear multiplication algorithms in ${\bbbf}_{q^n}$ over ${\bbbf}_q$.
\end{definition}

The complexity $\mu_{q,2}(n)$
will be denoted by $\mu_{q}(n)$, in accordance with the usual notation in the case of the product of two elements. 

\vspace{1em}

It will be interesting to relate the $k$-multilinear complexity of the multiplication to the minimal number $\nu_{q,k}(n)$ of 
bilinear multiplications in ${\bbbf}_q$ required to compute the product of
$k$ elements in the extension ${\bbbf}_{q^n}$. 

\begin{lemma}\label{Lemmetrivial}
\begin{equation}\label{Inega1}
\nu_{q,k}(1) \leq k-1,
\end{equation}

\begin{equation}\label{Inega2}
\nu_{q,k}(n) \leq (k-1) \times \mu_q(n),
\end{equation}

\begin{equation}\label{Inega3}
\mu_{q,k}(n) \leq \mu_{q,k}(mn) \leq \mu_{q,k}(m) \times \mu_{q^m,k}(n),
\end{equation}

\begin{equation}\label{Inega4}
\nu_{q,k}(n) \leq \nu_{q,k}(mn) \leq \nu_{q,k}(m) \times \nu_{q^m,k}(n).
\end{equation}

\end{lemma}
\begin{proof}
The two first inequalities are direct consequences
of the definitions. 
The inequalities (\ref{Inega3}) and (\ref{Inega4}) follow from the embedding of the field ${\bbbf}_{q^n}$ in the field ${\bbbf}_{q^{mn}}$. $\square$
\end{proof}

When for any $i$ the $k$ linear forms $a^*_{i,j}$ are all  the same linear form $a^*_{i}$,
namely

\begin{equation}\label{algo1}
t_{M_k}=\sum_{i=1}^{s} (a^*_{i})^{\otimes k} \otimes c_i
\end{equation}
and
$$\prod_{j=1}^k x_j=t_{M_k}(\otimes_{i=1}^k x_j)=\sum_{i=1}^{s}\left(\prod_{j=1}^ka^*_{i}(x_j)\right) c_i,$$
the decomposition is called a symmetric decomposition of $t_{M_k}$. 
As a consequence of \cite[Theorem 5]{bsho}, when $n>1$ such a symmetric decomposition exists if and only if $k \leq q$.
The symmetric $k$-multilinear complexity of the multiplication of $k$ elements in ${\bbbf}_{q^n}$ over ${\bbbf}_q$, 
denoted by $\mu^{Sym}_{q,k}(n)$,
is the minimum number of summands $s$ in the decomposition (\ref{algo1}). For $k=2$
this complexity is denoted by $\mu^{Sym}_{q}(n)$. From the definitions we get$\mu_{q,k}(n) \leq\mu^{Sym}_{q,k}(n)$.

\subsection{New results and organisation}

From a generalization of Chudnovsky-type algorithms to the k-multiplication, obtained by Randriambololona and Rousseau in \cite{rous1} (cf. also  \cite{rous2}), that we generalize to places of arbitrary degree, we obtain uniform upper bounds for the rank of the k-multiplication tensor in the finite fields (Theorem \ref{TheoremBoundsComp} in Section \ref{SectionBounds}). In this aim, we apply this type of algorithms to an explicit tower of Garcia-Stichtenoth \cite{gast2} and the corresponding descent tower described in Section \ref{towers}. Note that Randriambololona and Rousseau only obtain an asymptotic upper bound in $O(n)$ by using Shimura curves used by Shparlinski, Tsfasman and Vladut in \cite{shtsvl}.

\section{Theoretical construction of a multiplying algorithm}
\subsection{Notations}\label{ga}
Let $F/{\bbbf}_q$ be an algebraic function field over the finite field ${\bbbf}_q$
of genus $g$. We denote by $N_i(F/{\bbbf}_q)$ the number of places of degree $i$ of $F$ over
${\bbbf}_q$. If $D$ is a divisor, ${\mathcal L}(D)$ denotes the Riemann-Roch space
associated to $D$. Let $Q$ be a place of $F/{\bbbf}_q$. We denote by ${\mathcal O}_Q$ the valuation ring
of the place $Q$ and by $F_Q $ the residue class field ${\mathcal O}_Q/Q$ of the place $Q$ 
which is isomorphic to ${\bbbf}_{q^{{\rm deg} (Q)}}$ where ${\rm deg} (Q)$ is the degree of the place $Q$. 
Let us recall that for any $g \in {\mathcal O}_Q$, $g(Q)$ denotes the class of $g$  in ${\mathcal O}_Q/Q=F_Q$.  
Let us define the following Hadamard product in  ${\bbbf}_{q}^{N_1} \times {\bbbf}_{q^{2}}^{N_2} 
\times \cdots \times {\bbbf}_{q^{d}}^{N_d}$ 
where the $N_i$ denote integers $\geq 0$:

$$\bigodot_{i=1}^k(u_{1,1}^i,\cdots,u_{1,N_1}^i,\cdots, u_{d,1}^i,\cdots,u_{d,N_d}^i)= $$ 
$$\left(\prod_{i=1}^ku_{1,1}^i,\cdots,\prod_{i=1}^ku_{1,N_1}^i,\cdots, 
\prod_{i=1}^ku_{d,1}^i,\cdots,\prod_{i=1}^ku_{d,N_d}^i\right)
$$

\subsection{Algorithm}\label{ga1}
The following theorem 
generalizes the known results of the case $k=2$.

\begin{theorem} [Algorithm]\label{ChudGene}
Let
\begin{enumerate}
\item $q$ be a prime power and $k \geq 2$ be an integer,
\item $F/{\bbbf}_q$ be an algebraic function field,
\item $Q$ be a degree $n$ place of $F/{\bbbf}_q$,
\item ${\mathcal D}$ be a divisor of $F/{\bbbf}_q$,
\item ${\mathcal P}=\{P_{1,1},\cdots,P_{1,N_1}, \cdots, P_{d,1},\cdots,P_{d,N_d}\}$ 
be a set of $N=\sum_{i=1}^{d}N_i$ places of arbitrary degree where $P_{i,j}$ denotes a place of 
degree $i$ and $N_i$ a number of places of degree $i$.
\end{enumerate}
We suppose that $Q$ and all the places in $\mathcal P$ are not in the support of ${\mathcal D}$ and that:
\begin{enumerate}
\item the map
$$
Ev_Q:  \left \{
\begin{array}{ccl}
\Ld{} & \rightarrow & {\bbbf}_{q^n}\simeq F_Q\\
f & \longmapsto & f(Q)
\end{array} \right.
$$ 
is onto,
\item the map
$$
Ev_{\mathcal P} :  \left \{
\begin{array}{ccl}
\Ld{k} & \longrightarrow & {\bbbf}_{q}^{N_1} \times {\bbbf}_{q^{2}}^{N_2} \times \cdots \times {\bbbf}_{q^{d}}^{N_d} \\
f & \longmapsto & \left(\strut f(P_{1,1}),\cdots,f(P_{1,N_1}), \cdots, f(P_{d,1}),\cdots,f(P_{d,N_d})\right)
\end{array} \right.
$$
is injective
\end{enumerate}
Then, for any $k$ elements $x_1,....,x_k$ in ${\bbbf}_{q^n}$, we have 
$$m_k(x_1,...,x_k)= 
Ev_Q\left(Ev_{\mathcal P}^{-1}\left(\bigodot_{i=1}^{k} \left(Ev_{\mathcal P}
\left(\strut Ev_Q^{-1}(x_i)\right)\right)\right)\right),$$ 
and 
\begin{equation}\label{ineqcomp1}
\mu_{q,k}(n) \leq \sum_{i=1}^d N_i\mu_{q,k}(i) ,
\end{equation}

\begin{equation}\label{ineqcomp2}
\nu_{q,k}(n) \leq (k-1)\sum_{i=1}^d N_i\mu_{q}(i).
\end{equation}

\end{theorem}

\begin{proof}
For any $k$ elements $x_1,....,x_k$ in ${\bbbf}_{q^n}$, we have $k$ elements $g_1,\cdots, g_k$ in  $\Ld{}$ 
such that $Ev_Q(g_i)=g_i(Q)=x_i$ where $g_i(Q)$ denotes the class of $g_i$ 
in the residue class field $F_Q$.

Thus, 
\begin{equation}
m_k(x_1,...,x_k)=\prod_{i=1}^k x_i=\prod_{i=1}^kg_i(Q)=\left(\prod_{i=1}^kg_i\right)(Q).
\end{equation}
Moreover, since the divisor ${\mathcal D}$ has a positive dimension, without less of generality we 
can suppose that ${\mathcal D}$  is an effective divisor. 
Hence, $\Ld{}^k\subset \Ld{k}$ and $\prod_{i=1}^kg_i\in \Ld{k}$. 
Now, we have to compute the function $h=\prod_{i=1}^kg_i$ and this is done by an interpolation process 
{\it via} the evaluation map $Ev_{\mathcal P}$.
More precisely, $$h=\prod_{i=1}^kg_i=Ev_{\mathcal P}^{-1}\left(\bigodot_{i=1}^{k}Ev_{\mathcal P}(g_i)\right),$$
namely$$h = Ev_{\mathcal P}^{-1}\left(\bigodot_{i=1}^{k} 
\left(Ev_{\mathcal P}\left(\strut Ev_Q^{-1}(x_i)\right)\right)\right).$$ Moreover,  
$\Ld{k}\subset {\mathcal O}_Q$ for any integer $k$, thus, we can apply the map $Ev_Q$ 
over the element $h$ which corresponds to 
the product of $k$ elements $x_1, ...,x_k$ which gives the first inequality. Inequality (\ref{ineqcomp2}) follows from Inequality 
(\ref{Inega2}) in Lemma \ref{Lemmetrivial} and we are done. $\square$
\end{proof}

\begin{theorem} \label{theoprinc}
Let $q$ be a prime power and let $n$ be an integer $>1$. 
Let $F/{\bbbf}_q$ be an algebraic function field of genus $g$.
Let ${\mathcal P}_i$ be a set of places of degree $i$ in $F/{\bbbf}_q$ and $N_i$ 
the cardinality of ${\mathcal P}_i$. We denote
${\mathcal P}=\bigcup_{i=1}^r{\mathcal P}_i$.
Let us suppose that there is a place of degree $n$ and a non-special divisor of degree $g-1$.
If there is an integer $r\geq 1$ such that:
\begin{equation}\label{maincondition}
\sum_{i=1}^r iN_i > kn+kg-k
\end{equation}
then $$\mu_{q,k} (n) \leq \sum_{i=1}^r 	N_i \,\mu_{q,k}(i).$$
Let $r_0, r'_0$ such that 
$$\frac{\mu_{q,k}(r_0)}{r_0}=\sup_{1\leq i \leq r}\frac{\mu_{q,k}(i)}{i} \hbox{ and }  %Max_{1\leq i \leq r}\left(\strut\mu_{q,k}(i)\right)
\frac{\mu_{q}(r'_0)}{r'_0}=\sup_{1\leq i \leq r}\frac{\mu_{q}(i)}{i}.$$
Then,
\begin{equation}\label{borne1} 
\mu_{q,k} (n) \leq \left(kn+kg-k+r\right) \,\frac{\mu_{q,k}(r_0)}{r_0}. 
\end{equation}

\begin{equation}\label{borne2} 
\nu_{q,k} (n) \leq (k-1) \left(kn+kg-k+r\right) \,\frac{\mu_{q}(r'_0)}{r'_0}. 
\end{equation}

\end{theorem}

\begin{proof}
Let $\mathcal R$ be a non-special divisor  
of degree $g-1$.
Then we choose a divisor ${\mathcal D}_1$ such that ${\mathcal D}_1={\mathcal R}+Q$. 
Let $[{\mathcal D}_1]$ be the class of ${\mathcal D}_1$, then by 
\cite{deur}, Lecture 14, Lemme 1, $[{\mathcal D}_1]$ contains a divisor $\mathcal D$ 
defined over ${\bbbf}_q$ such that $ord_{P}{\mathcal D}=0$ for all places $P \in {\mathcal P}$ 
and $ord_{Q}{\mathcal D}=0$. Since $ord_Q{\mathcal D}=0$, 
${\mathcal L}({\mathcal D})$ is contained in the valuation ring $O_Q$ of $Q$. 
Hence $Ev_Q$ is a restriction of the residue class mapping, 
and defines an ${\bbbf}_q$-algebra homomorphism. 
The kernel of $Ev_Q$ is ${\mathcal L}({\mathcal D}-Q)$. But ${\mathcal D}-Q$ is non-special divisor of degree $g-1$, 
then $l({\mathcal D}-Q)=deg({\mathcal D}-Q)-g+1=0$ and $Ev_Q$ is injective with $deg~{\mathcal D}=n+g-1$.
Moreover, if $\mathcal K$ is a canonical divisor, we have 
$l({\mathcal D})=l({\mathcal K}-{\mathcal D})+n$ by the Riemann-Roch theorem.
Hence, $l({\mathcal D}) \geq n$ and as $Ev_Q$ is injective, 
we obtain $l({\mathcal D})=n$. We conclude that $Ev_Q$ is an isomorphism. 
The map $Ev_{\mathcal P}$ is well-defined because 
${\mathcal L}(k{\mathcal D})$ is contained in the valuation ring 
of every place of ${\mathcal P}$. The kernel of $Ev_{\mathcal P}$ is 
${\mathcal L}\left(k{\mathcal D} - \sum_{P\in {\mathcal P}}P\right )$ 
which is trivial because $\sum_{i=1}^r iN_i > kn+kg-k$. Therefore $Ev_{\mathcal P}$ 
is injective and $l(k{\mathcal D})=kn+(k-1)g-k+1$.

Let us remark that we can suppose that $kn+kg-k< \sum_{i=1}^r i N_i \leq kn+kg-k+r$.
Then, by Inequality (\ref{ineqcomp1}) in Theorem \ref{ChudGene}, 
$$\mu_{q,k} (n) \leq \sum_{i=1}^r N_i\mu_{q,k} (i) = \sum_{i=1}^r i N_i\frac{\mu_{q,k}(i)}{i}\leq \frac{\mu_{q,k}(r_0)}{r_0}\sum_{i=1}^r i N_i,$$
which gives Inequality (\ref{borne1}). Inequality (\ref{borne2}) follows from Inequality (\ref{ineqcomp2}) in Theorem \ref{ChudGene} in the same way and the proof is complete. $\square$
\end{proof}

\begin{remark}
A place of degree $n$ and a non-special divisor of degree $g-1$ 
exist if elementary numerical conditions are satisfied. More precisely, if 
\begin{equation}\label{existences} 
2g+1\leq q^{\frac{n-1}{2}}(q^{\frac{1}{2}}-1),
\end{equation} 
then there exists a place of degree $n$ by \cite[Corollary ]{stic2}. Moreover, if $q\geq 4$ or $N_1\geq g+1$ 
then there exits a non-special divisor of degree $g-1$ by \cite{balb} (cf. also \cite{bariro}).
\end{remark}

\begin{remark}
Note that it is an open problem to know if the bilinear complexity 
is increasing or not with respect to the extension degree, when $n\geq 2$. 
\end{remark}

\section{On the existence of these algorithms}

\subsection{Strategy of construction}\label{strat}

In this section, we present our strategy of construction of these algorithms which is well adapted to an asymptotical study.
More precisely, this strategy consists in fixing the definition field ${\bbbf}_q$, the integer $n$ and the parameter $k$ of the algorithm \ref{ChudGene}. 
Then, we suppose the existence of family of algebraic function fields over ${\bbbf}_q$ of the genus $g$ growing to the infinity.

\vspace{1em}

The main condition (\ref{maincondition}) of Theorem \ref{theoprinc} supposes that we can find algebraic function fields 
having good properties. In particular, it is sufficient to have a family of function fields having 
sufficiently places with a certain degree $r$. In this aim, we focalize on sequences of algebraic functions fields with increasing genus 
attaining the Drinfeld-Vladut bound of order $r$ (cf. \cite[Definition 1.3]{baro4}). 
Hence, let us find the minimal integer $r$ for such a family, in order to  
Condition (\ref{maincondition}) is satisfied. 

\vspace{1em}

Let ${\mathcal G}/{\bbbf}_q= (G_i/{\bbbf}_q)_i$ be a sequence of algebraic function fields $G_i$ over ${\bbbf}_q$ of genus $g_i$ attaining 
the Drinfeld-Vladut Bound of order $r$, namely $$\lim_{i\rightarrow +\infty}\frac{N_r(G_i)}{g_i}=\frac{1}{r}(q^{\frac{r}{2}}-1).$$
The divisors $D$ are such that $\deg (D)=n+g-1$, thus $\deg kD =kn+kg-k$. To determine the minimal value of the integer $r$, we can suppose 
without less of generality that $rN_r>kn+kg-k$, which implies asymptotically that $q^{\frac{r}{2}}-1>k$, namely $r>2\log_q(k+1)$.

\subsection{Towers of algebraic function fields}\label{towers}

In this section, we present a sequence of algebraic function fields defined over ${\bbbf}_{q}$ from the Garcia-Stichtenoth 
tower constructed in \cite{gast2}.  This tower is suitable with respect to the strategy defined in Section \ref{strat}.

\vspace{1em}
 
Let  us consider a finite field ${\bbbf}_{l^2}$ where $l$ is a prime power such that ${\bbbf}_{l^2}$ is an extension field of ${\bbbf}_q$. 
We consider the Garcia-Stichtenoth's elementary abelian tower $\mathcal{F}$ over ${\bbbf}_{l^2}$ constructed in \cite{gast2} 
and defined by the sequence $\mathcal{F}=(F_0, F_1,\cdots, F_i,\cdots )$ where 
$$
F_0 := {\bbbf}_{l^2}(x_0)
$$
is the rational function field over ${\bbbf}_{l^2}$, and for any $i \geq 0$, $F_{i+1} := F_i(x_{i+1})$ with $x_{i+1}$ satisfying the following equation:  
$$
x^l_{i+1}+x_{i+1} = \frac{x^l_i}{x_i^{l-1} + 1}.
$$
\noindent Let us denote by $g_i=g(F_i)$ the genus of $F_i$ in $\mathcal{F}/{\bbbf}_{l^2}$ and recall the following formul\ae:
\begin{equation}\label{genregs}
g_i = \left \{ \begin{array}{ll}
		(l^{\frac{i+1}{2}}-1)^2  & \mbox{for odd } i, \\
		(l^{\frac{i}{2}}-1)(l^{\frac{i+2}{2}}-1) & \mbox{for even } i. 
		\end{array} \right .
\end{equation}

\noindent Thus, as in \cite{bapi2}, according to these formul\ae, it is straightforward that the genus of any step of the tower satisfies:
\begin{equation}\label{bornegenregs}
(l^{\frac{i}{2}}-1)(l^\frac{i+1}{2}-1) < g_i < (l^{\frac{i+2}{2}}-1)(l^\frac{i+1}{2}-1).
\end{equation}
Moreover, a tighter upper bound will be useful and can be obtained by expanding expressions in (\ref{genregs}):
\begin{equation}\label{bornesupgenregs}
g(\mathbf{F}_i) \leq l^{i+1}-2l^\frac{i+1}{2}+1.
\end{equation}

Then we can consider as in \cite{balbro} the descent tower  $\mathcal{G}/{\bbbf}_{q}$ defined over ${\bbbf}_{q}$ given by the sequence:
$$G_0 \subset G_1 \subset \cdots \subset G_i \subset \cdots$$
defined over the constant field ${\bbbf}_{q}$ and related to the tower $\mathcal{F}$ by: 
$F_i={\bbbf}_{l^2}\otimes_{{\bbbf}_{q}} G_i$ for all i. 

\vspace{1em}

Let us recall the known results concerning the number of places of degree one of 
the tower~$\mathcal{F}/{\bbbf}_{l^2}$, established in \cite{gast2}. 

\begin{proposition}\label{genus}
The number of places of degree one of $F_i/{\bbbf}_{l^2}$ is :
$$
N_1(\mathbf{F}_i/{\bbbf}_{l^2}) = \left\{ 
	\begin{array}{ll}
	 l^i(l^2-l)+2l^2 & \mbox{if the caracteristic is even},\\
	 l^i(l^2-l)+2l & \mbox{if the caracteristic is odd}.
\end{array}\right.
$$
\end{proposition}

Now, we are interested by supplementary properties concerning the descent tower $\mathcal{G}/{\bbbf}_{q}$ in the context of Theorem \ref{theoprinc}, which we are going to apply to this tower. More precisely, we need that the steps of the tower $\mathcal{G}/{\bbbf}_{q}$ verify the mandatory properties of the existence of a place of degree $n$ and of sufficient number of places of certain degrees. 
In this aim, we introduce the notion of the action domain of an algebraic function field.

\begin{definition} 
Let us define the following quantities:
\begin{enumerate}
\item $M_i=l^i(l^2-l)$ where $l=q^{\frac{r}{2}}$,
\item $\Delta_{q,k,i}=M_i-kg_i+k$,
\item $\Theta_{q,k,i}=\{n \in {\bbbn} \mid \Delta_{q,k,i}>kn\},$ 
\item $R_{q,k,i}=\sup \Theta_{q,k,i},$
\item the set $\phi_{q,i}=\{n \in {\bbbn}\mid 2g_i+1\leq l^{\frac{n-1}{r}}(l^{\frac{1}{r}}-1)\},$
\item $\Gamma_{q,i}=\inf \phi_{q,i}.$
\item $I_{q,k,i}=\Theta_{q,k,i}\cap \phi_{q,i}=[\Gamma_{q,i},R_{q,k,i}]$
\end{enumerate}
\end{definition}

\begin{lemma}
Let $q$ be a prime power and $k \geq 2$ be an integer and $r$ the smallest even integer 
$>  2\log_q(k+1)$.
Then $(\Delta_{q,k,i})_{i\in {\bbbn}}$ is an increasing sequence such that 
$\lim_{i\rightarrow \infty}\Delta_{q,k,i}=+\infty$.
\end{lemma}

\begin{proof}
By Inequality (\ref{bornesupgenregs}) we get:
$$\Delta_{q,k,i}=M_i-kg_i+k\geq q^{\frac{ri}{2}}\left( q^{r}-q^{\frac{r}{2}} \right) 
-k \left(q^{\frac{r(i+1)}{2}}-2q^{r\frac{i+1}{4}} \right).$$ 
Hence:
$$\Delta_{q,k,i}\geq q^{\frac{ri}{2}}(q^{r}-q^{\frac{r}{2}})-kq^{\frac{r(i+1)}{2}}=
q^{\frac{r(i+1)}{2}}(q^{\frac{r}{2}}-1-k).$$ 
As  $r > 2\log_q(k+1)$, 
we are done. $\square$
\end{proof}

\begin{proposition}
Let $q$ be a prime power and $k \geq 2$ be an integer and $r$ the smallest 
even integer $>  2\log_q(k+1)$.
Then for any $i\geq 1$, the action domain $I_{q,k,i}$ of $F_i/{\bbbf}_q$ 
is not empty.
\end{proposition}

\begin{proof}
Let us compute bounds on the two values $\Gamma_{q,i}$ and $R_{q,k,i}$.
\begin{enumerate}
 \item Bound on $\Gamma_{q,i}$. The set $\phi_{q,i}$ contains the set of integers $n$ such that
 \begin{equation}\label{contrainte2}
 2\left(l^{i+1}-2l^{\frac{i+1}{2}}+1\right)+1 \leq l^{\frac{n-1}{r}}\left(l^{\frac{1}{r}}-1\right).
 \end{equation}
 If 
 $$n \geq r(i+1)+1 +\log_{l^{\frac{1}{r}}}(2)-\log_{l^{\frac{1}{r}}}(l^{\frac{1}{r}}-1),$$
 Condition (\ref{contrainte2}) is verified,
 then the integer $n$ is in $\phi_{q,i}$. We conclude that
 $$\Gamma_{q,i} \leq r(i+1)+3.$$
 
 \item Bound on $R_{q,k,i}$. The set $\Theta_{q,k,i}$ contains the set of integers $n$ such that
 \begin{equation}\label{contrainte1}
  q^{\frac{r(i+1)}{2}}(q^{\frac{r}{2}}-1-k)\geq kn.
 \end{equation}
As $r > 2\log_q(k+1)$, the integer $q^{\frac{r}{2}}-1-k$ is $\geq 1$.
Then the condition 
$$n \leq \frac{1}{k}q^{\frac{r(i+1)}{2}}$$
implies Condition (\ref{contrainte1}).
As $q^{\frac{r}{2}}>k+1$ the condition
$$n \leq (k+1)^i$$
implies $n \in \Theta_{q,k,i}$.
Then
$$R_{q,k,i}\geq (k+1)^i.$$
 \item Conclusion. For $i \geq 1$, the function $(k+1)^i-r(i+1)-3$ is an increasing
 function of $i$. For $i=1$ we have  $[2r+3,k+1] \subset [\Gamma_{q,1},R_{q,k,1}]$,
 and more generally $[r(i+1)+3,(k+1)^i]\subset [\Gamma_{q,i},R_{q,k,i}]$.
\end{enumerate}

$\square$

\end{proof}

\begin{proposition}

Let $q$ be a prime power and $k \geq 2$ be an integer and $r$ the smallest 
even integer $>  2\log_q(k+1)$. Then:
$$[2r+3,+\infty[ \subset \bigcup_{i \geq 1} I_{q,k,i}.$$ 
\end{proposition}

\begin{proof}
For any $i$, we have $[r(i+1)+3,(k+1)^i]\subset I_{q,k,i}$ and $\bigcup_{i \geq 1} I_{q,k,i}=[2r+3,+\infty[$ and the proof is complete.
$\square$

\end{proof}

Hence, for any integer $n\geq 2r+3$, there exists an action domain $I_{q,k,i}$ such that $n\in I_{q,k,i}$. 

\section{Uniform upper bounds}\label{SectionBounds}

In this section, we study the complexity of the family of algorithms \ref{ChudGene} constructed with a 
tower $\mathcal{G}/{\bbbf}_{q}$ defined over ${\bbbf}_{q}$ and defined in Section \ref{towers}, such that 
this tower attains the Drinfeld-Vladut bound of order $r$ where $r$ is the smallest even integer $> 2\log_q(k+1)$. 
Moreover, for any integer $i$, we have $\sum_{j\mid r}jN_j(G_i)=N_1(F_i)$. 

\begin{theorem}\label{TheoremBoundsComp}
Let $q$ be a prime power and $k \geq 2$ be an integer and $r$ the smallest even integer $>  2\log_q(k+1)$. 
Let $r_0, r'_0$ such that
$$\frac{\mu_{q,k}(r_0)}{r_0}=\sup_{1\leq i \leq r}\frac{\mu_{q,k}(i)}{i} 
\hbox{ and } \frac{\mu_{q}(r'_0)}{r'_0}=\sup_{1\leq i \leq r}\frac{\mu_{q}(i)}{i}.$$
Then for any integer $n$, we have: 

\begin{equation}\label{Formcomp1}
\mu_{q,k} (n) \leq \frac{k(kq^{\frac{r}{2}}+1)n-k+1}{r_0} \mu_{q,k}(r_0)
\end{equation}
and consequently 

\begin{equation}\label{Formcomp2}
\nu_{q,k} (n) \leq  \frac{k(k-1)(kq^{\frac{r}{2}}+1)n-(k-1)^2}{r'_0} \mu_{q}(r'_0).
\end{equation}

\end{theorem}

\begin{proof}
We apply Theorem \ref{theoprinc} to the tower descent tower $\mathcal{G}/{\bbbf}_{q}$ by taking for any $n>1$ the small step $F_i$ satisfying the assumptions of Theorem \ref{theoprinc}.
Let us set $M_i=l^i(l^2-l)$ where $l=q^{\frac{r}{2}}$.
For any integer $k \geq 2$ and any integer $n$, let $i$ be the smallest integer such that $M_i> kn+kg_i-k$, then $kn<M_i-kg_i+k$
and $kn\geq M_{i-1}-kg_{i-1}+k$. So, we obtain 
$$kn\geq q^{\frac{r(i-1)}{2}}\left( q^{r}-q^{\frac{r}{2}} \right) 
-k \left(q^{\frac{ri}{2}}-2q^{\frac{ri}{4}} \right)$$ 
by 
Formula (\ref{bornesupgenregs}). Hence, $kn\geq q^{\frac{ri}{2}} \left( q^{\frac{r}{2}}-k-1\right)$, and so:
$$q^{\frac{ri}{2}}\leq \frac{kn}{q^{\frac{r}{2}}-k-1}.$$
But $q^{\frac{r}{2}}-k-1\geq 1$, hence: 
\begin{equation}\label{imin}
i\leq \frac{2}{r}\log_q(kn).
\end{equation}
Now, by Bound (\ref{borne1}) in Theorem \ref{theoprinc}, we obtain:
$$\mu_{q,k} (n) \leq \frac{kn+kg_i-k+r}{r_0}  \,\mu_{q,k}(r_0)$$ 
which gives by Formulae (\ref{bornesupgenregs}) and (\ref{imin}): 
$$\mu_{q,k} (n) \leq \frac{kn+kq^{\frac{r(i+1)}{2}}-k+r}{r_0}  \,\mu_{q,k}(r_0)$$ 
with 
$$q^{\frac{r(i+1)}{2}}\leq q^{\log_q(kn)+\frac{r}{2}}=q^{\frac{r}{2}}kn,$$
which gives the first inequality. Inequality (\ref{Formcomp2}) follows from Inequality (\ref{borne2}) in Theorem \ref{theoprinc}.
$\square$
\end{proof}

\begin{corollary}
Let $q$ be a prime power and $k \geq 2$ be an integer and $r$ the smallest even integer 
$> 2\log_q(k+1)$. Let $r_0, r'_0$ such that
$$\frac{\mu_{q,k}(r_0)}{r_0}=\sup_{1\leq i \leq r}\frac{\mu_{q,k}(i)}{i} \hbox{ and } %Max_{1\leq i \leq r}\left(\strut\mu_{q,k}(i)\right)
\frac{\mu_{q}(r'_0)}{r'_0}=\sup_{1\leq i \leq r}\frac{\mu_{q}(i)}{i}.$$

Then for any integer $n$, we have: 

\begin{equation} \label{Ineq1Coro}
\mu_{q,k} (n) \leq \frac{k(kq^{\frac{r}{2}}+1)}{r_0} \mu_{q,k}(r_0)\,n\leq \frac{k(k(k+1)q+1)}{r_0} \mu_{q,k}(r_0)\,n
\end{equation}

\begin{equation} \label{Ineq2Coro}
\nu_{q,k} (n) \leq  \frac{k(k-1)(kq^{\frac{r}{2}}+1)}{r'_0} \mu_{q}(r'_0)\,n\leq\frac{k(k-1)(k(k+1)q+1)}{r'_0} \mu_{q}(r'_0)\,n .
\end{equation}

\end{corollary}

\begin{proof}
The first inequality of (\ref{Ineq1Coro}) follows immediately from Inequality (\ref{Formcomp1}) in Theorem \ref{TheoremBoundsComp}.
As $r$ the smallest even integer $> 2\log_q(k+1)$, we have $r-2<2\log_q(k+1)$ and so $q^{r/2-1}<k+1$.
Hence, $q^{r/2}<q(k+1)$ which gives the first assertion. Assertion (\ref{Ineq2Coro}) is obtained with the same way.
$\square$
\end{proof}

\end{document}